\documentclass[11pt]{article}
\usepackage{epsfig}
\usepackage{psfrag}
\usepackage{MnSymbol}

\title{A robust spectral method for finding lumpings and meta stable states of non-reversible Markov chains\thanks{This work was funded by PACE (Programmable Artificial Cell Evolution), a European Integrated Project in the EU FP6-IST-FET Complex Systems Initiative, by EMBIO (Emergent Organisation in Complex Biomolecular Systems), a European Project in the EU FP6-NEST Initiative, and by MORPHEX (Morphogenesis and gene regulatory networks in plants and animals: a complex systems modelling approach), a European Project in the EU FP6-STREP Initiative.}}

\author{Martin Nilsson Jacobi\thanks{Department of Energy and Environment, Chalmers University of Technology, Gothenburg, Sweden, (mjacobi@chalmers.se)} }

\begin{document}

\maketitle
\begin{abstract}
A spectral method for identifying lumping in large Markov chains is presented. Identification of meta stable states is treated as a special case. The method is based on spectral analysis of a self-adjoint matrix that is a function of the original transition matrix. It is demonstrated that the technique is more robust than existing methods when applied to noisy non-reversible Markov chains.
\end{abstract}

%\begin{keywords}
 %Markov chain, stochastic matrix, metastable states, lumping, aggregation, modularity, block diagonal dominance, block stochastic 
%\end{keywords}

%\begin{AMS}
%15A18, 15A51, 60J10, 65F15 
%\end{AMS}

\section{Introduction}

The structural dynamics of large biomolecules can often be accurately described as a Markov transition process.  Frequently, the dynamics display separation of time scales where aggregated conformational states are evolving at much slower rate than the detailed molecular dynamics does. The problem of identifying the conformational states from the detailed Markov transition matrix has received recent interest~\cite{deuflhard00identification,Deuflhard05,Christian2008,friszache}. The technically similar problem of identifying modularity and community structure on complex networks has also been discussed extensively, e.g.~\cite{fiedler,pothen,newman}.

Identification of meta stable states is a special case of a more general reduction called (approximate) lumping. Lumping of a Markov chain means that the state space is partitioned into equivalence classes of states called macro states. A coarse grained process is defined by the transitions between the macro states. If the coarse grained process is Markovian, i.e. exhibits no memory, we call the reduction a lumping. A partition into meta stable states is an example of an approximate lumping in the following sense. In the limit of complete stability, i.e. when there are no transitions between the macro states, then the macro states define a degenerate case of exact  lumping. More generally, the Markov property is fulfilled on the aggregated level if the relaxation process within a meta stable state is fast and mixing so that the memory of exactly how the meta stable state was entered is lost before the transition to a new meta stable state occurs. In this sense aggregation into meta stable states can be viewed as an approximate lumping. Aside from separation of time scales, there are other generic situations when a Markov process is expected to be lumpable. For example when a particle  interacts with many other particles, a ``heat bath'', the dynamics of the single particle can be described as a Brownian motion. Technically the transition matrix of a lumpable Markov chain can be rearranged into a block-stochastic structure, see Fig.~\ref{fig:block_stoch} and definition (\ref{eq:block_stoch}). Markov chains with metastable states can be permuted into a block-diagonal structure (Fig.~\ref{fig:block_dominant}), which is a special case of a block-stochastic matrix.

The most successful methods for identifying  meta stable states and modules in networks are based on the level structure of the eigenvectors whose corresponding eigenvalues are clustered  close to the Perron-Frobenius eigenvalue. The technique introduced in this paper is closely related to these spectral method first introduced by Fidler in the 70's, at that point as a method for graph partitioning~\cite{fiedler}. Fidler noted that the second eigenvector of the graph Laplacian shows tightly connected communities of nodes that are connected to the  other communities by relatively few edges, or low algebraic connectivity. Later the method was used in connection to the classic graph coloring problem~\cite{aspvall}.  In the same paper the idea of using the sign structure of the $k$ first eigenvectors to partition a graph into $k$ aggregates of nodes was introduced. The same idea was later applied to identify meta stable states in Markov chains~\cite{deuflhard00identification}. For these spectral  methods to be stable, the eigenvalue problem must be symmetric with respect to some scalar product. This means that the Markov process must be reversible, or that the network is assumed to be effectively non-directional. A notable exception in presented based symmetrization using the stationary distribution was presented in~\cite{Froyland}. Another exception is a recent method for Markov chains based on singular value decomposition of the Markov transition matrix~\cite{friszache}. However, the SVD-based method is not appropriate for identifying lumpings of Markov chains since the singular vectors typically do not have a level structure (or relevant sign structure) in the theoretical limit of exact lumpability, see for example the transition matrix defined in~(\ref{eq:simpleEx}).

In this paper we present a new robust spectral method for identifying possible lumpings of non-reversible Markov chains. Instead of using the spectrum of the transition matrix directly, we define a self-adjoint ``invariance matrix'' whose kernel relates to the eigenvectors that define the meta stable states, or more generally the lumps of the Markov chain. Since the invariance matrix is self-adjoint by construction, the usual assumption of reversibility can be lifted. We demonstrate the method of both Markov chains with meta stable states and more general block-stochastic structure, and compare the performance to other methods reported in the literature, e.g. the methods presented in~\cite{Froyland} and~\cite{friszache}.

\section{Lumping of Markov chains}
\label{sec:lumping}

Consider a Markov process $x_{t+1} = x_t P$. The $N \times N$ transition matrix $P$ is a row stochastic matrix, i.e. $\sum _j P_{ij} = 1$ $\forall i$. A lumping is defined as a partition of the states space $\Sigma$ into $K$ equivalence classes of states $L_k$ such that $L_k \cap L_l = \emptyset$ and $\cup _k L_k = \Sigma$~\cite{Rogers}. A necessary and sufficient condition for a partition to be a lumping is~\cite{Kemeny99}
\begin{equation}
	\sum _{j \in L_l} P _{ij} \;\;\;\; \mbox{constant for all $i$ in an aggregate, $i \in L_k$.}
	\label{eq:lumping}
\end{equation}
If a Markov chain allows for a non-trivial lumping we call it lumpable. A simple example of a lumpable transition matrix is
\begin{equation}
	P = \frac{1}{4} \left( \begin{array}{ccc} 3 & 0 & 1 \\ 1 & 2 & 1 \\ 0 & 2 & 2 \end{array} \right) ,
		\label{eq:simpleEx}
	\end{equation}
which, aside from the trivial lumping defined by all states aggregated into one macro-state, allows the non-trivial lumping $\{ \{ 1 , 2 \} , 3 \}$, i.e. state $1$ and $2$ lumped into one macro-state. 

The condition in (\ref{eq:lumping}) also immediately defines the transition matrix for the aggregated dynamics
\begin{equation}
	\widetilde{P} _{kl} = \sum _{j \in L _l} P_{ij} \;\;\;\; \mbox{$i \in L_k$},
	\label{eq:pt}
\end{equation}
since all states $i \in L_k$ give the same result. In practice Eq.~\ref{eq:lumping} is usually not fulfilled exactly. For example, if a transition matrix can be written as
\begin{equation}
	P = (1-\epsilon) A + \epsilon B ,
\end{equation}
where $A$ is a transition matrix that fulfills the lumpability condition (\ref{eq:lumping}) and $B$ is some arbitrary transition matrix. Then, if $\epsilon$ is small, we say that $P$ is approximately lumpable. Note that the aggregated transition probabilities in (\ref{eq:lumping}) are in this case approximately constant with deviations of $\mathcal{O} ( \epsilon )$. The reduced dynamics must in this case be approximated e.g. using a weighted average for the transitions between the aggregated states
\begin{equation}
	\widetilde{P} _{kl} = \frac{1}{\sum _{j \in L_l} v_j } \sum _{i \in L_k} \sum _{j \in L _l} v_j P_{ij} ,
	\label{eq:effective_dyn}
\end{equation}
where $v_j$ is the stationary distribution. Using the weighted average is natural since it gives the same reduced transition matrix as we find if we estimate the aggregated transition probabilities directly from a stationary time series.

A partition can be represented by a matrix $\Pi$ defined as $\Pi _{ik} = 1$ if $i \in L_k$ and $\Pi _{ik} =0$ otherwise. Eq. (\ref{eq:lumping}) can be reformulated as
\begin{equation}
	P \Pi = \Pi \widetilde{P} ,
	\label{eq:commuting}
\end{equation}
which, if written out explicitly in terms of the elements, implies that the column space of $\Pi$ spans a right-invariant subspace of $P$. Assuming that $P$ is diagonalizable, the invariant subspace is spanned by a set of right eigenvectors of $P$, and due to the $0$ or $1$ structure of $\Pi$ the elements in these eigenvectors must be constant over the aggregates. To be more precise, a lumping with $K$ aggregates exists if and only if there are exactly $K$ right eigenvectors of $P$ with elements that are constant over the aggregates, see~\cite{Shi01arandom,jacobi} for details. As an example, the transition matrix defined in (\ref{eq:simpleEx}) allows for the lumping $\{ \{ 1, 2 \}, 3 \}$ as indicated by the two first elements in the right eigenvectors $( 1 ,1 ,1 ) ^T $ and $( -1 , -1 , 2) ^T$ being constant.

It should be noted that there exist other types of aggregation of states where the aim is to preserve (for example) the structure of the equilibrium distribution. A prominent example of this is renormalization of lattice spin systems. However, in this paper we focus on lumping that respect the {\em dynamics} of the process. In this case the Markov property is the central constraint, i.e. the mutual information between the past and the future given the present should be zero on both the micro (a prerequisite for the procedure) and macro level (the lumping condition). This leads to the strong conditions on the aggregation seen in Eq.~\ref{eq:lumping} and Eq.~\ref{eq:commuting}. For a more detailed discussion on how memory appears on the coarse grained level if the lumping criterion is not fulfilled, see~\cite{jacobi}.

The principle idea behind spectral methods for identifying lumping or meta stable states, as well as modularity in networks, is to search for (right) eigenvectors whose elements are constant over the aggregates, i.e. eigenvectors with a level structure, see e.g.~\cite{jacobi}   for details. If the transition matrix is symmetric under some scalar product the eigenvectors are orthogonal and it is easy to show that  the constant level sets must have different sign structure~\cite{aspvall} (the sign structure of a vector is defined by mapping negative elements to $-1$ and positive elements to $+1$). The sign structure is often used as a lumping criterion rather than the constant levels since this is expected to be a more numerically stable~\cite{aspvall,deuflhard00identification}. However, a more recent study has shown that the sign structure is more sensitive to noise than the constant level structure over the aggregates, an observation that lead the authors to introduce an algorithm based on the simplex structure of the almost constant level sets~\cite{Deuflhard05}.

For the detection of metastable states or modularity, the eigenvectors of interest are those corresponding to eigenvalues close to the Perron-Frobenius eigenvalue, since these eigenvalues are related to the slow dynamics. In the case of general lumping the eigenvectors involved are not necessarily distinguished by their appearance in the spectrum. However, as we discuss in Section~\ref{sec:block_stoch_matrix}, for large transition matrices the eigenvectors involved in lumping tend to be separated from the rest of the spectrum by being located further from the origin in the complex plane than the rest of the spectrum, but not necessarily by being closer to $1$.

As a complement to the spectral methods, the commutation relation (\ref{eq:commuting}) can be used directly to identify lumping of Markov chains. Start by making a random assignment of the $N$  states to $K$ aggregates, and construct the corresponding $\Pi$ matrix. Given the $\Pi$ matrix, the reduced transition matrix $\widetilde{P}$, defined in Eq. \ref{eq:pt} can be derived by simply ignoring that the row elements are not constant within the aggregates and use the average defined in (\ref{eq:effective_dyn}).
The left hand side in (\ref{eq:commuting}), $P \Pi$, defines a $K$ dimensional vector for each of the $N$ states. If the lumping is correct then all states in an aggregate $k$ should have identical $K$ dimensional vectors, and they should be equal to the $k$th row of $\widetilde{P}$. If $\Pi$ is not a lumping we can try to improve it by assigning state $i$ to aggregate $k$ where $k = \mbox{argmin} _l \| ( P \Pi ) _i - \widetilde{P} _l \|$. The result is a new aggregation with a new $\Pi$ matrix. The process can be iterated until convergence. A similar method was introduced by Lafon and Lee~\cite{Lafon}. As pointed out in \cite{weinan} it is similar in structure to the $K$-means clustering algorithm. It should be noted that this direct clustering technique  only works if the aggregated dynamics has long relaxation time, i.e. there is a spectral gap supporting the lumping, see~\cite{Lafon}   for details.
%This can be realized by noting that the algorithm works by iteratively applying $P$ so that the dominant modes that define the aggregation become fixated. 
The performance of the algorithm is shown in comparison with the method introduced in this paper in Fig.~\ref{fig:results_block_diag} and~\ref{fig:results_block_stoch}.

\section{A robust method for identifying lumping}

We now present the main idea of this paper. We would like to find invariant vectors containing invariant level sets. For moderately sized (unperturbed) transition matrices or for time-reversible Markov chains, the eigenvectors can be used to detect lumping. If the Markov chain is not reversible, calculation of both eigenvalues and eigenvectors is numerically unstable, since, for example, the transition matrix may contain non-trivial Jordan blocks~\cite{jacobi}. This is the motivation for the new method.  Start by noting that, if a normalized vector $u$ is approximately right invariant under $P$, then there must exist a $\lambda$ such that
\begin{equation}
	\| (P - \lambda I ) u \| _2 \ll  1 ,
	\label{baru}
\end{equation}
where $I$ denotes the identity matrix.	The square of the $2$-norm on the left hand side in (\ref{baru}) is not sensitive to small changes in the elements of $P$, whereas if $P$ is non-symmetric the eigenvalues and eigenvectors can be ill conditioned. Obviously, if $u$ is an eigenvector and $\lambda$ the corresponding eigenvalue, then (\ref{baru}) is zero, reflecting the fact that the eigenvector is exactly invariant. The $2$-norm of  (\ref{baru}) is given by
\begin{equation}
	u ^{\dagger} Q ( \lambda ) u ,
	\label{min}
\end{equation}
where $u ^{\dagger}$ denotes the conjugated transpose of the vector $u$. The ``invariance matrix'' $Q$ is defined as
\begin{equation}
	Q ( \lambda ) = P ^{\dagger} P - \lambda ^* P - \lambda P ^{\dagger} + | \lambda | ^2 I ,
	\label{def_Q}
\end{equation}
(note that $Q$ is typically not a stochastic matrix). Regardless of the properties of $P$, $Q ( \lambda) $ is by construction a self-adjoint matrix and diagonalization is numerically stable. If $\lambda$ is an eigenvalue of $P$, then $Q (\lambda )$ is positive semi-definite with a zero eigenvalue corresponding to the eigenvector of $P$ with eigenvalue $\lambda$. If $\lambda$ is not an eigenvalue, then $Q (\lambda )$ is positive definite. For a given $\lambda$, (\ref{baru}) is minimized by $u$ being the eigenvector of $Q (\lambda )$ corresponding to the smallest eigenvalue of $Q ( \lambda )$, or in the case of degeneracy a linear combination of the eigenvectors of the smallest eigenvalue. 

%It is important to note that due to the intrinsic instability of the eigenvalues of a non-symmetric matrix, (\ref{baru}) can be small even when $\lambda$ is not close to any of the actual eigenvalues of $P$. In many situations approximate lumpings are revealed by eigenvectors of $Q ( \lambda )$ where $\lambda$ is chosen relatively far from any of the eigenvalues of $P$. In these cases the $\lambda$'s that reveal the approximate lumping are eigenvalues of another, unperturbed,  transition matrix that exactly accepts the lumping in question.

\section{Meta stable states}

\begin{figure} 
	\psfrag{err}[][]{$\;$\scriptsize $\tau$}
	\psfrag{eps}[][]{\scriptsize{$\epsilon$}}
	\begin{center}
		\epsfig{file=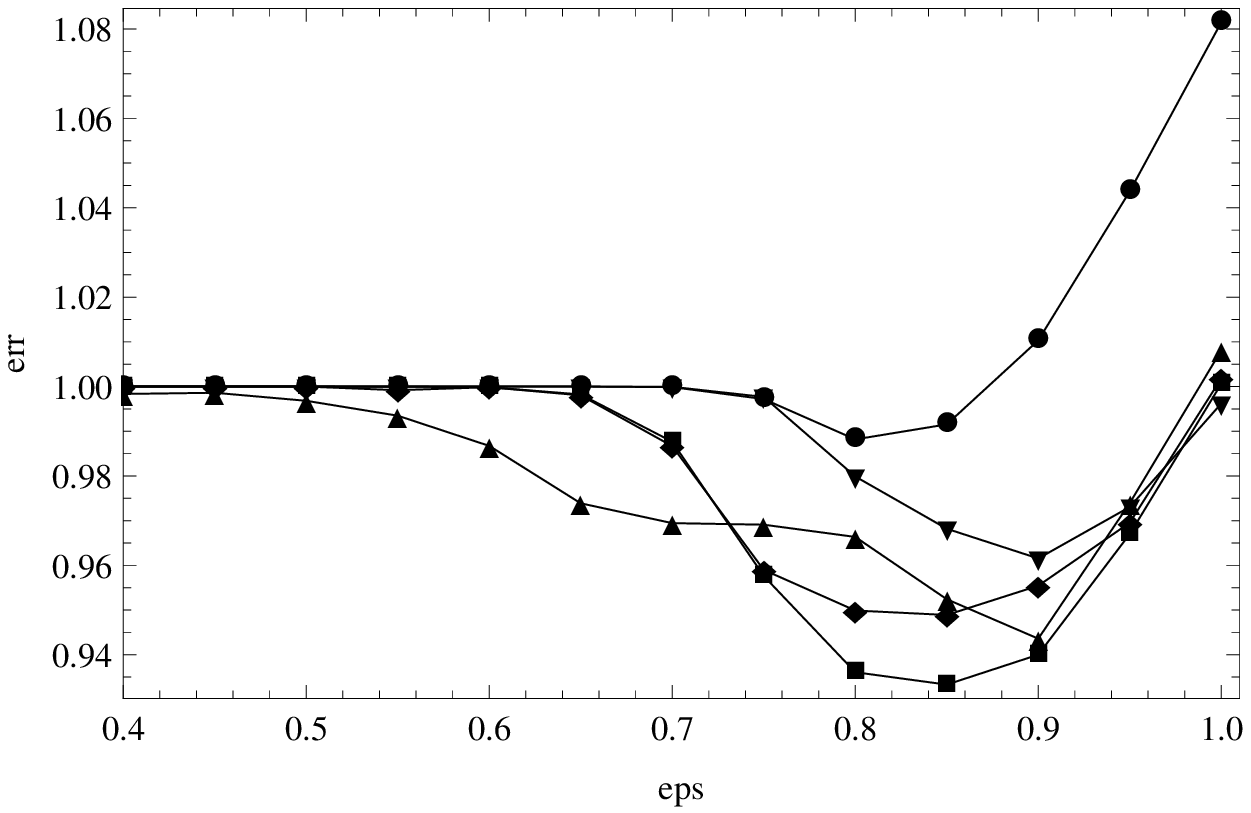,width=6cm} \hspace{0.5cm}
 \epsfig{file=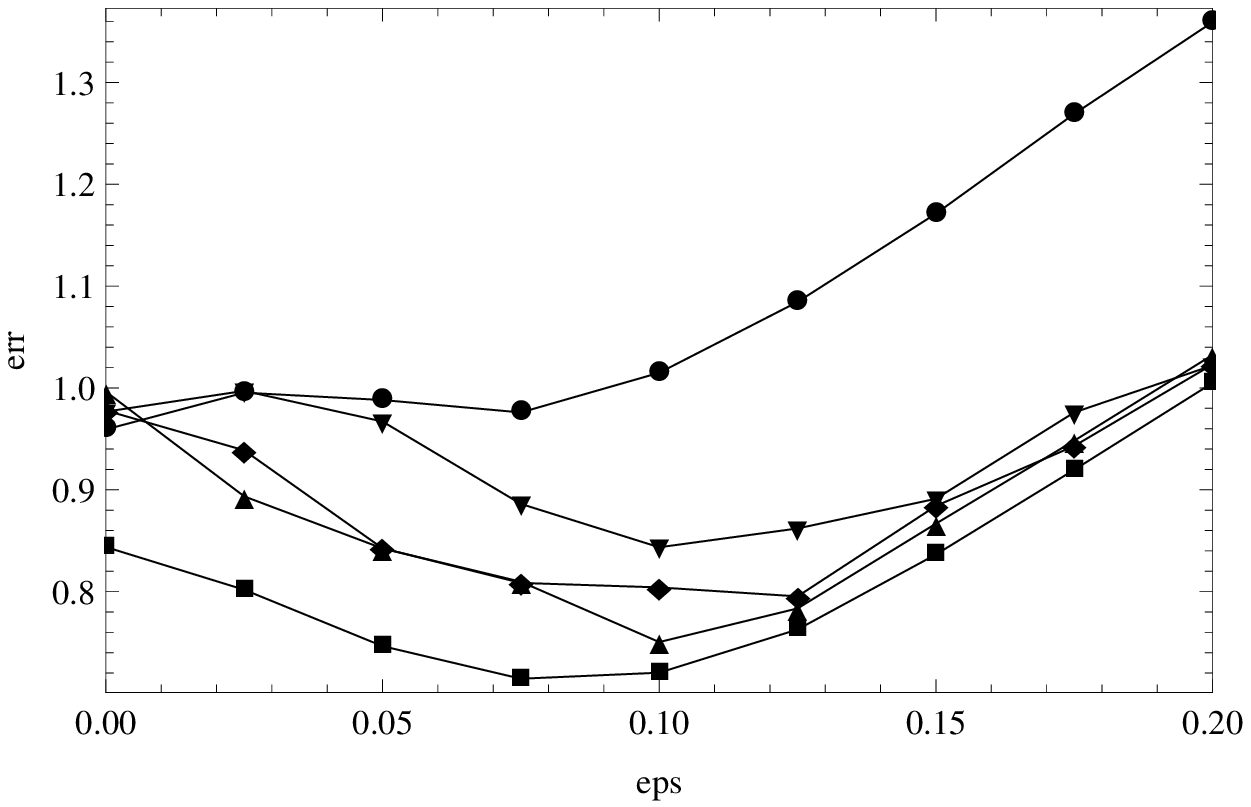,width=6cm}
\end{center}
\caption{The average result of identifying the meta stable states of a dominant block diagonal transition matrices defined in (\ref{eq:block_diag}) is shown on the left, and for matrices defined in (\ref{eq:block_diag2}) with $\delta = 0.2$ on the right. The average waiting time for transitions between meta stable states (normalized against the result for the partitioning used to generate the matrices) for $100$ test matrices of size $200 \times 200$ for each $\epsilon$ value is used as a measure of the quality of the results. The result for the $Q$ method is displayed as $\bullet$, the result for using eigenvectors of $P$ as $\filledmedsquare$, the singular vectors from SVD \cite{friszache} as $\filleddiamond$, results from the clustering method suggested in~\cite{Lafon} (see Sec.~\ref{sec:lumping} for details) as $\filledmedtriangledown$, and the method suggested in~\cite{Froyland} as $\filledmedtriangleup$.}
\label{fig:results_block_diag}
\end{figure}

Detecting meta stable states is an especially simple, but also especially interesting, case of general (approximate) lumping. The meta stable states are characterized by their long relaxation time, and hence their dynamics is associated with eigenvalues close to $1$. The right eigenvectors involved in the aggregation have corresponding eigenvalues closer to $1$ than the rest of the spectrum. As a consequence, meta stable states can be identified by the approximate constant level structure of the eigenvectors of 

\begin{equation}
	Q (1) =  P ^{\dagger} P -  P -  P ^{\dagger} + I 
	\label{Q_1}
\end{equation}
with eigenvalues close to $0$. As a consequence, the small eigenvalues of (\ref{Q_1}) include the eigenvectors needed in the aggregation. It should be noted that the actual eigenvalues of $P$ associated with the meta stable states need not be close to $1$ for the sub-dominant eigenvectors of $Q(1)$ to reveal the meta stable states, see Fig.~\ref{fig:eig_vec_block_diag} and \ref{fig:spectrum_block_diag}. Since $Q$ is self-adjoint the eigenvectors are orthogonal. The eigenvectors are approximately constant over the aggregates, and orthogonality can then only be achieved if each aggregate has a unique sign structure in the eigenvectors. This observation was used by Aspvall and Gilbert~\cite{aspvall} and Deuflhard and coworkers~\cite{deuflhard00identification,Deuflhard05}, but in these cases under the condition of symmetry of the adjacency matrix or reversibility of the transition matrix respectively. Using the $Q$ matrix there is no need to make assumptions on $P$. It is straight forward to apply the same sign structure criterion to the eigenvectors of $Q$, but empirical tests have shown that in our case the following simple approach is relatively robust (see Fig.~\ref{fig:results_block_diag}): 

\begin{enumerate}
	\item Find the eigenvectors $\{ u_i \} _{i=1} ^K$ corresponding to the small eigenvalues of $Q(1)$ in (\ref{Q_1}).
	\item For each state $j = 1 , \dots , N$ form a $K$ dimensional vector $u_{\bullet j} = ( u_{1j} , u_{2j} , \dots , u_{Kj})$ of its corresponding elements in the $K$ eigenvectors. 
	\item Use a standard clustering algorithm (e.g. K-mean) to cluster the states with respects to the $u_{\bullet j}$ vectors. Note that we expect the level structure in the eigenvectors to be relatively stable to perturbations ($\mathcal{O} ( \epsilon ^2 )$), as pointed out in~\cite{Deuflhard05}.
\end{enumerate}

\begin{figure} 
\begin{center}
	\epsfig{file=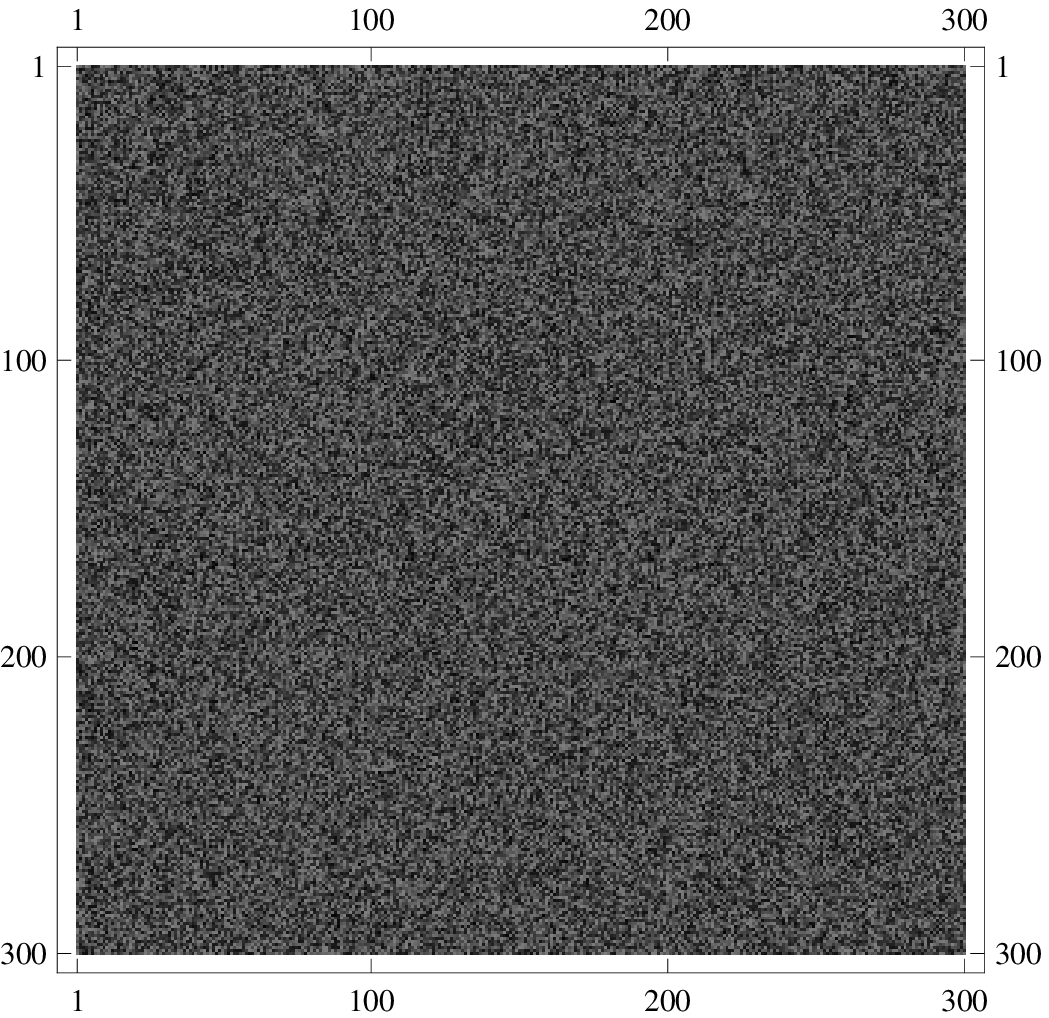,width=6cm}  
  \epsfig{file=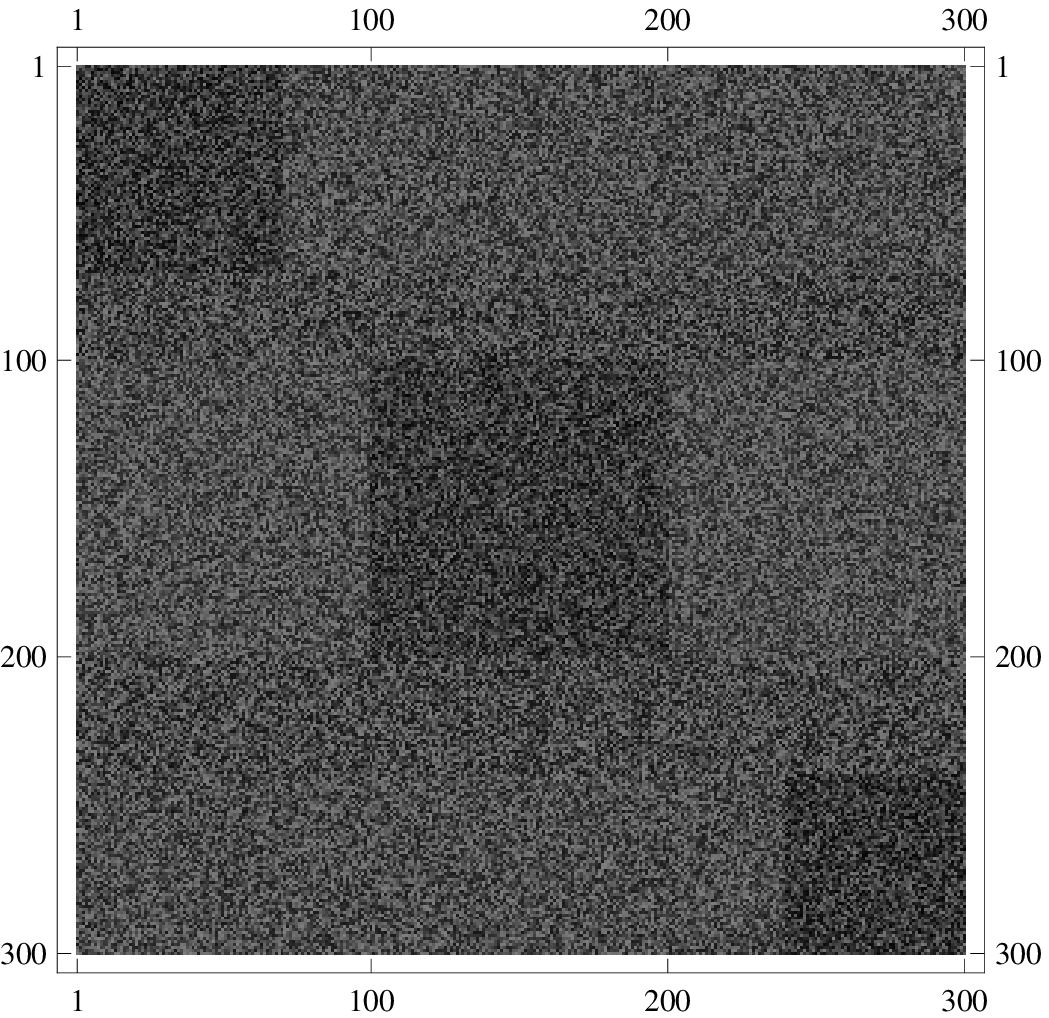,width=6cm} 
\end{center}
\caption{A weakly block dominant transition matrix, constructed as (\ref{eq:block_diag}) with $\epsilon = 0.7$, is shown. The time scale separation is not very pronounced, see the spectrum in Fig.~\ref{fig:spectrum_block_diag}. To the left with random permutation and to the right after sorting the matrix according to the aggregation revealed in the clusters of the eigenvectors of the $Q$ matrix shown in Fig.~\ref{fig:eig_vec_block_diag}.}
\label{fig:block_dominant}
\end{figure}

\begin{figure}
\psfrag{u1}[][]{$\;\;\;q_1$}
\psfrag{u2}[][]{$\;q_2$}
\begin{center}
	\epsfig{file=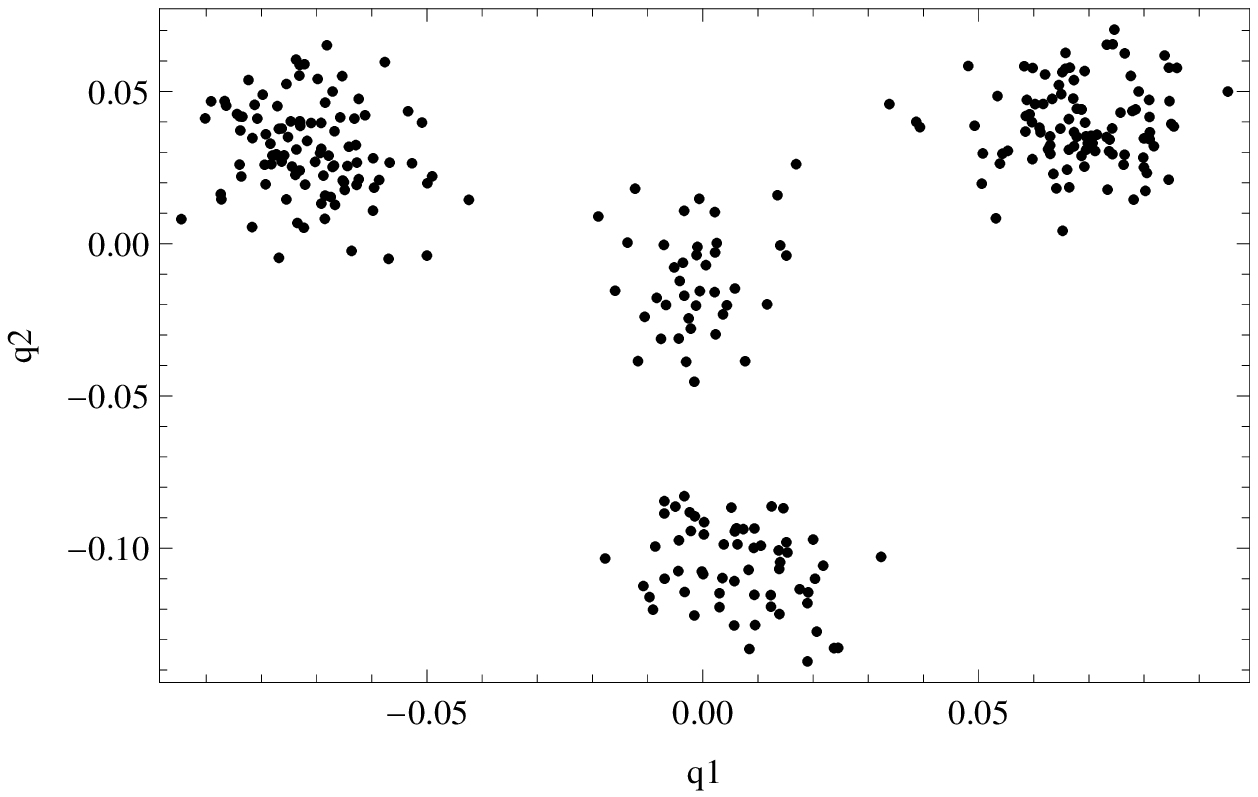,width=6cm} \hspace{0.3cm} 
	\psfrag{u1}[][]{$\;\;\;u_1$}
  \psfrag{u2}[][]{$\; u_2$}
	\epsfig{file=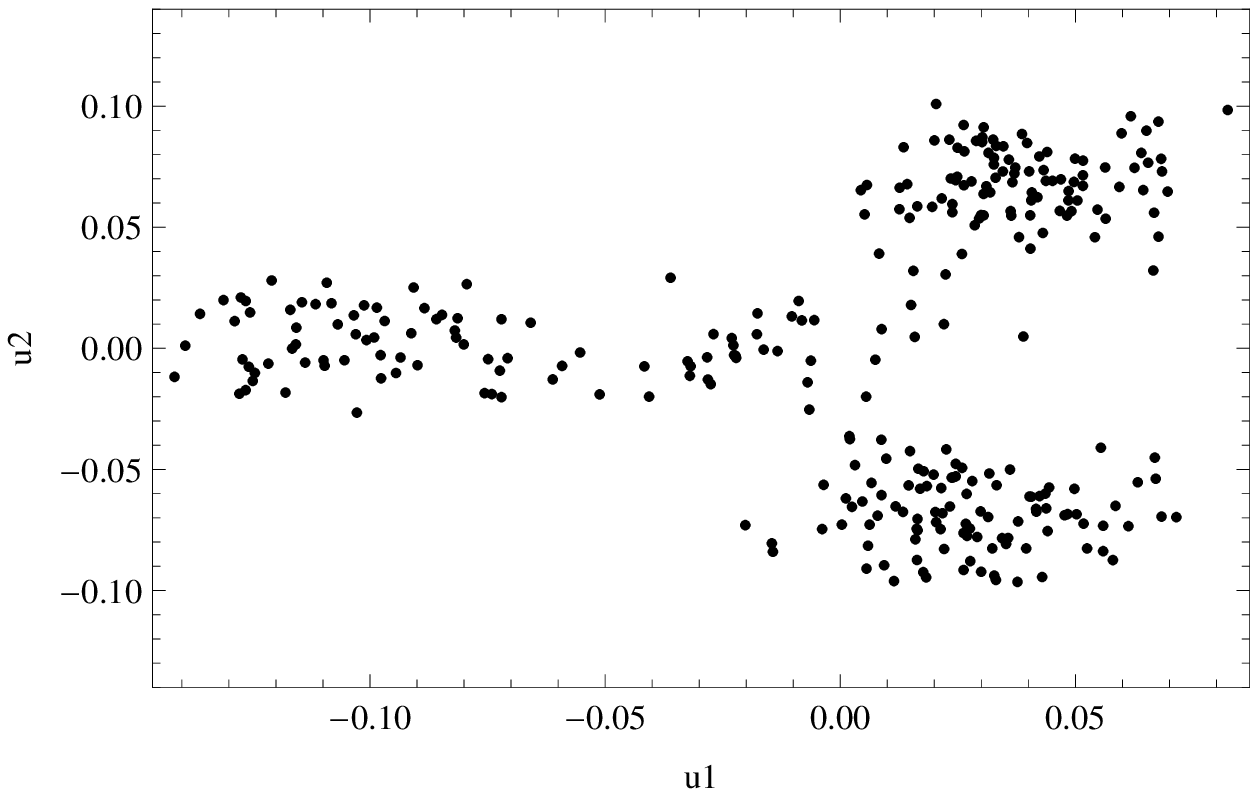,width=6cm} 
\end{center}
\caption{The figure shows the clustering of the elements in the second and third smallest, respective largest, eigenvectors of $Q(1)$ (to the left) respectively $P$ (to the right) of the matrix shown in Fig.~\ref{fig:block_dominant}. Note that the clusters are more distinct in the eigenvectors of $Q (1)$ shown on the left.}
\label{fig:eig_vec_block_diag}
\end{figure}

\begin{figure} 
\begin{center}
	\psfrag{Re}[][]{$\;$\scriptsize{Re($\lambda$)}}
	\psfrag{value}[][]{\scriptsize{$\lambda$}}
	\psfrag{Im}[][]{$\;$\scriptsize{Im($\lambda$)}}
\psfrag{index}[][]{\scriptsize{index}}
	\epsfig{file=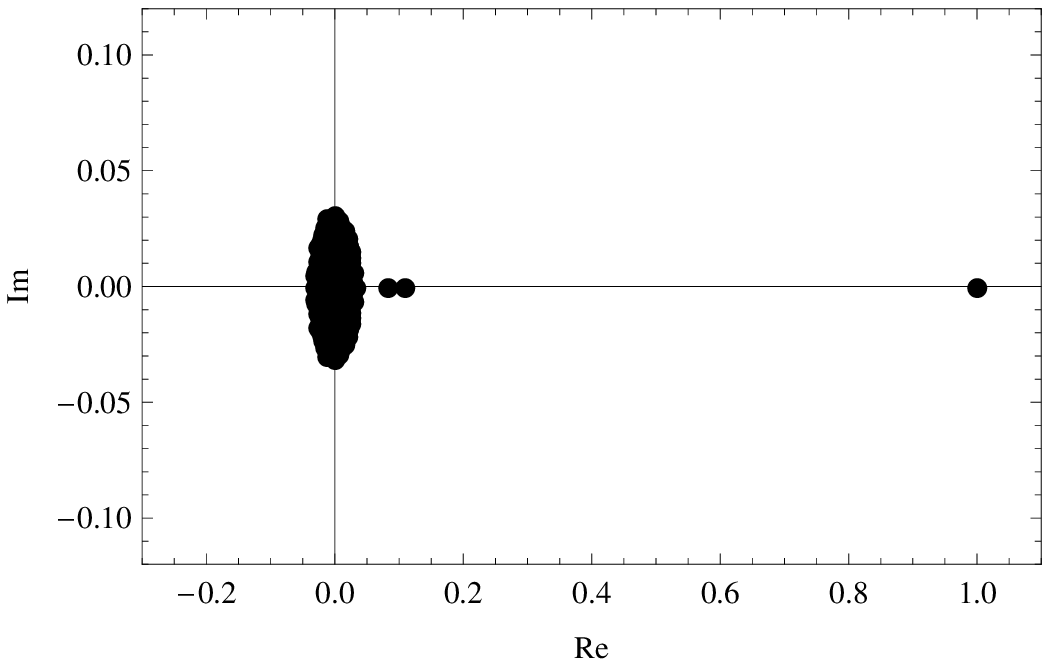,width=6cm}\hspace{0.5cm} 
  \epsfig{file=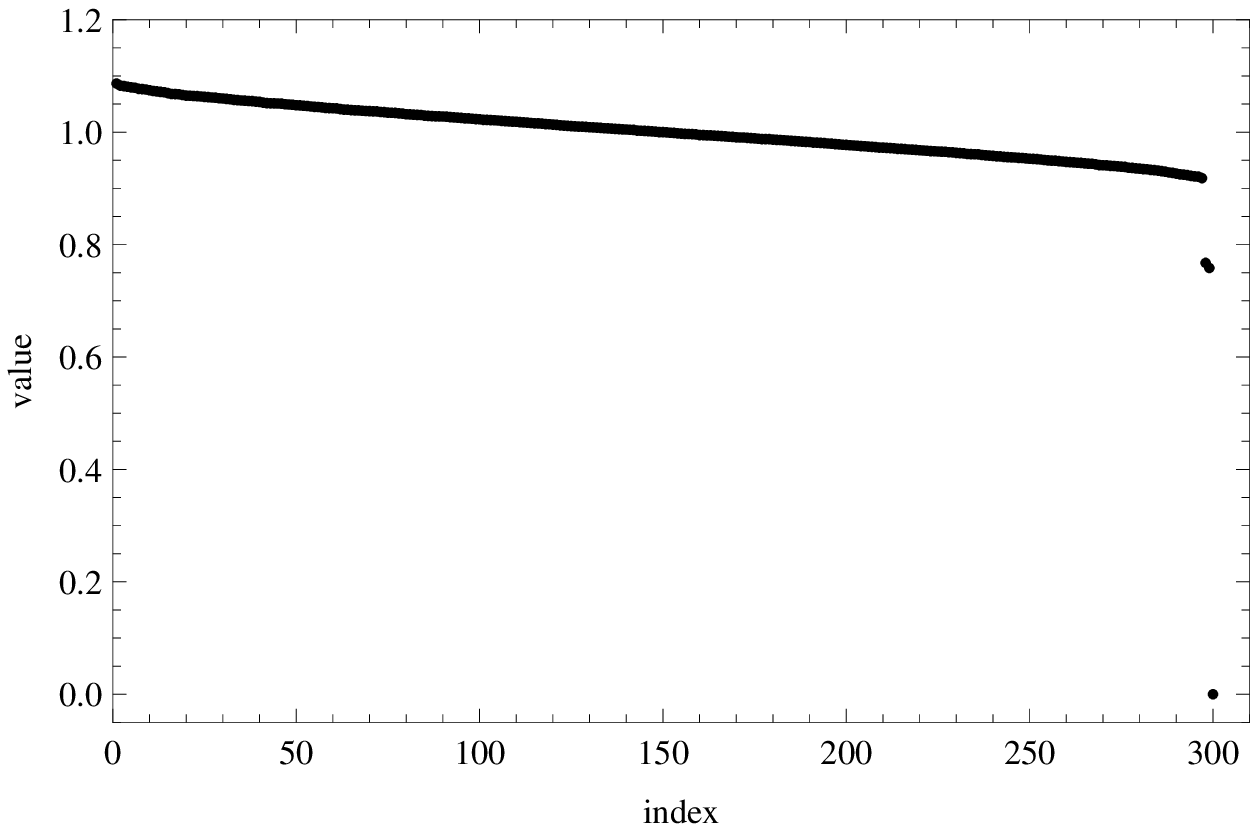,width=6cm} 
\end{center}
\caption{The spectrum of the transition matrix $P$ in Fig.~\ref{fig:block_dominant} to the left and of the corresponding $Q(1)$ matrix on the right. The aggregation into meta stable states is associated with the dominant eigenvalues of $P$, i.e. the Perron-Frobenius eigenvalue and the two eigenvalues close to $0.1$, or alternatively with the three smallest eigenvalues of $Q(1)$ (to the far right in the figure). Note that even though the dominant eigenvalues of $P$ are not clustered close to $1$, the spectrum and eigenvectors of $Q(1)$ show the meta stable states more clearly than the eigenvectors of $P$, see Fig.~\ref{fig:eig_vec_block_diag}.}
\label{fig:spectrum_block_diag}
\end{figure}

To test the method we generate two classes of matrices. The first class is on the form
\begin{equation}
	P = (1 - \epsilon ) B + \epsilon A ,
	\label{eq:block_diag}
\end{equation}
where $B$ is a block diagonal transition matrix with $3-5$ blocks and transition probabilities within the blocks chosen uniformly in the interval $[ 0 , 1]$ and then normalized. The matrix $A$ is a transition matrix with no block diagonal structure, generated in the same way as the blocks in $B$. The parameter $\epsilon$ sets the level of perturbation of $P$ from being block diagonal. Fig.~\ref{fig:block_dominant}--\ref{fig:spectrum_block_diag} show an example of a transition matrix of this type with $\epsilon = 0.7$ and the corresponding spectrum and clustering of the elements in the  dominant eigenvectors of $P$ and the sub-dominant eigenvectors of $Q$. 

It can be argued that matrices of the type in (\ref{eq:block_diag}) are unlikely to appear in practical applications. Instead of the smooth average modulation that the decrease transition probabilities between blocks in (\ref{eq:block_diag}), a more binary modulation often occurs in practice, i.e. many transition probabilities are zero. In this situation meta stable states occur as a consequence of a higher probability of having transitions within (rather than between) meta stable states. Contrasting the construction in (\ref{eq:block_diag}) this produces a sparse transition matrix. We construct the second class of matrices according to
\begin{equation}
	P ^* _{ij} ( \epsilon , \delta ) = \chi _{ij} (\epsilon , \delta ) B _{ij} ,
	\label{eq:block_diag2}
\end{equation}
where $B$, as before, is a matrix with random entries chosen uniformly in the interval $[ 0 , 1]$. In the matrix $\chi ( \epsilon , \delta )$, the entries are binary chosen so that $\chi _{ij} = 1$ with probability $\delta$ if $i$ and $j$ are in the same block,  and $\chi _{ij} = 1$ with probability $\epsilon$ if $i$ and $j$ are {\em not} in the same block, otherwise $\chi _{ij} = 0$. Thus $\delta$ controls the overall probability of transitions within the blocks and $\epsilon$ controls the transitions between the blocks, with $\delta \geq \epsilon$. The two extreme points are $\epsilon = 0$ which produce a completely block diagonal matrix, and $\epsilon = \delta$ which gives a matrix without any block diagonal structure. The rows in the matrix $P^*$ is normalized to produce a stochastic matrix. The procedure described can produce states with no outgoing transitions, i.e. that are ill defined as transition matrices. If this happens we generate a new matrix. 

We tested the performance of the $Q$ method and compared it to the following existing techniques: results from the eigenvectors of $P$, the right and left singular vectors from an SVD as suggested in~\cite{friszache}, the clustering method presented in~\cite{Lafon}, and the spectral method described in~\cite{Froyland}. As test cases we used  the two classes of matrices described above and measured how stable the meta stable states produced by the different methods where, i.e. the average waiting time between jumps between meta stable states. The result is shown in Fig.~\ref{fig:results_block_diag}. For each value of $\epsilon$ shown, the average switching time is measured for $100$ matrices of size $200 \times 200$ in respective class. The time $\tau$ is scaled so that the ``correct partitioning'' used to generate the matrices have waiting time $1$. The results indicate that the $Q$ method is more robust against perturbations than previously reported methods. It is especially interesting to note that for high $\epsilon$ values, i.e. when the original block diagonal structure is almost lost, the $Q$ method still produce aggregations that are more stable than random partitions. Non of the other methods are capable of finding these very weak meta stable states. 

\section{Block stochastic matrix}
\label{sec:block_stoch_matrix}
\begin{figure} 
\begin{center}
	\epsfig{file=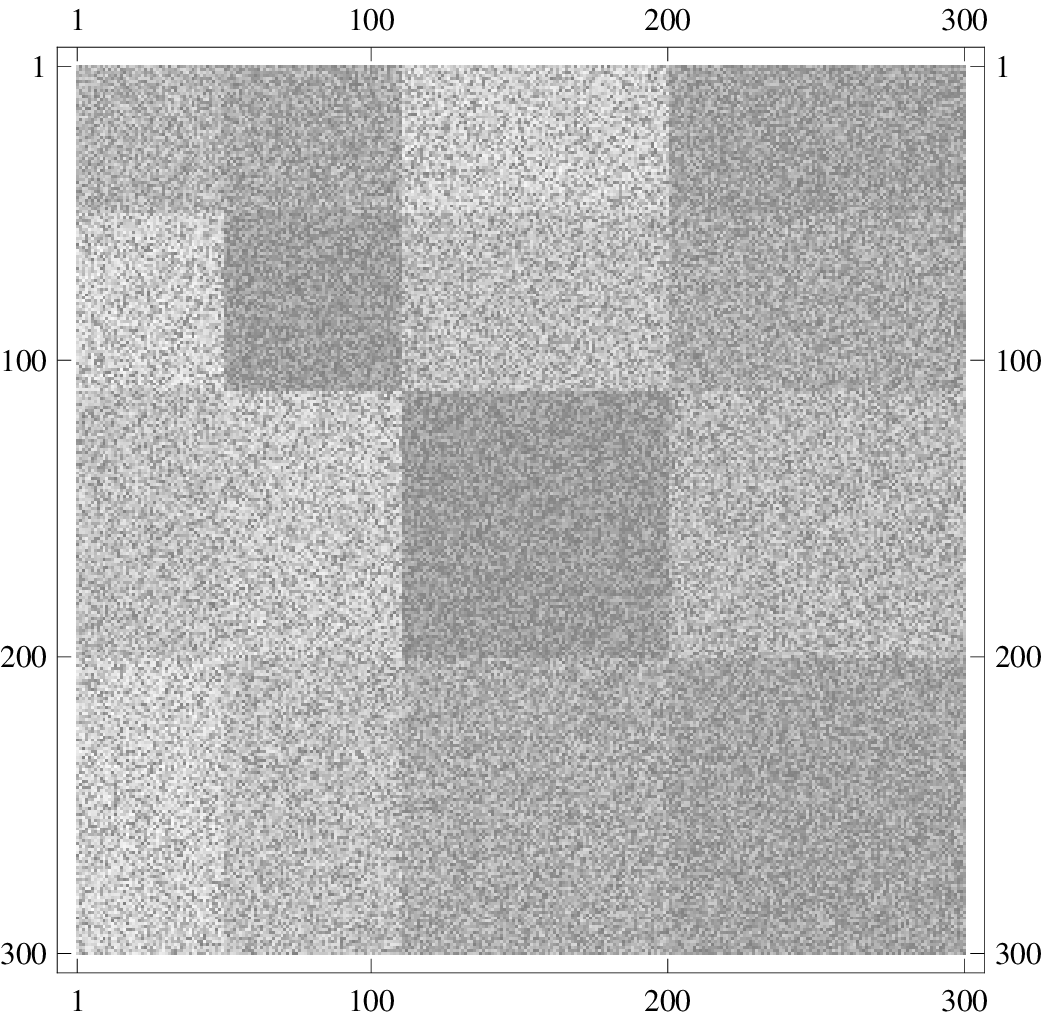,width=6cm}  
  \epsfig{file=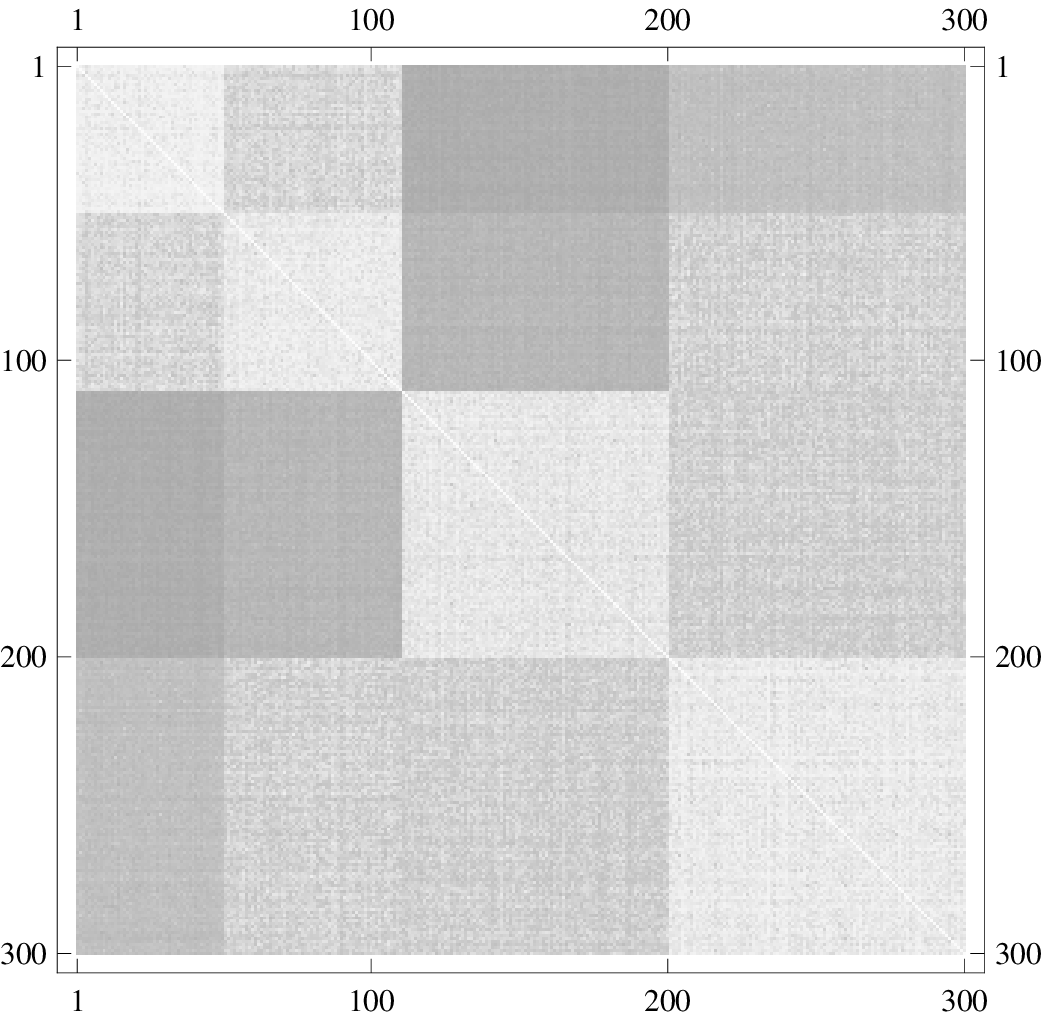,width=6cm} 
\end{center}
\caption{A block stochastic transition matrix $P$, constructed as (\ref{eq:block_stoch}) with $\epsilon = 0.5$, is shown to the left and $P P ^{^\mathrm{T}}$ to the right. The SVD method from~\cite{friszache} is based on the right matrix and it is clear that the two matrices share the same block stochastic structure. However the numerical tests show that $Q$ method based on the left matrix is more stable to perturbations than the SVD method based on the eigenvectors of the right matrix. The spectrum of $P$ is shown in  Fig.~\ref{fig:spectrum_block_stoch}. }
\label{fig:block_stoch}
\end{figure}

As discussed earlier, matrices with dominant block diagonal structure are special cases of lumpable Markov processes. The more general structure of lumpable transition matrices is shown in Fig.~\ref{fig:block_stoch}. A block-stochastic matrix is a matrix on the form
\begin{equation}
	P = \left( \begin{array}{cccc} \widetilde{P} _{11} a_{11} & \widetilde{P} _{12} a_{12} & \cdots & \widetilde{P} _{1m} a_{1m} \\
		\vdots & \vdots & \ddots & \vdots \\
		\widetilde{P} _{k1} a_{k1} & \widetilde{P}_{k2} a_{k2} & \cdots & \widetilde{P}_{mm} a_{mm} \end{array}
		\right) ,
		\label{eq:block_stoch}
	\end{equation}
where $\widetilde{P}$ is the $m \times m$ transition matrix of the reduced dynamics and each of the $a_{ij}$ is a transition matrix in itself. Naturally, $a_{ij}$ and $a_{ji}$ must for a fixed $i$ have the same dimensions for $j=1 , \dots , m$.

\begin{figure} 
\psfrag{Re}[][]{$\;$\scriptsize{$\;\;\;\;$Re($\lambda$)}}
	\psfrag{value}[][]{\scriptsize{$\lambda$}}
	\psfrag{Im}[][]{$\;$\scriptsize{Im($\lambda$)}}
	\begin{center}
  \epsfig{file=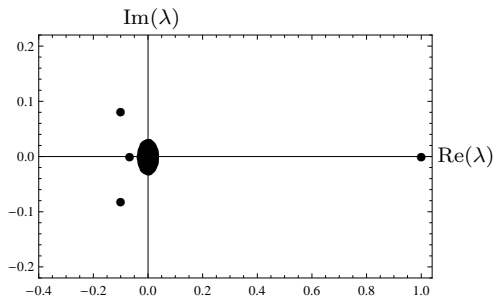,width=6cm} 
\end{center}
\caption{The figure shows the spectrum of the block stochastic matrices shown in Fig.~\ref{fig:block_stoch}. The eigenvalues associated with the eigenvectors that are involved in the lumping process are separated from the rest of the spectrum, but typically not close to the Perron-Frobenius eigenvalue. It should be noted that though the separation in the spectrum typically appears for large block stochastic matrices, this is a statistical effect and not necessary for the transition matrix to be lumpable, see the main text for further discussion on this. }
\label{fig:spectrum_block_stoch}
\end{figure}

The spectra of large block stochastic matrices tend to separate out the eigenvalues associated with the lumping. The separation is however different from the one occurring in block diagonally dominant transition matrices. Instead of clustering around the Perron-Frobenius eigenvalue, the reducing eigenvalues of a block stochastic matrix separate by larger distance to the origin in the complex plane, see Fig.~\ref{fig:spectrum_block_stoch}. The reason behind the separation is also different from the block diagonal case, and only appears as a statistic effect for large random transition matrices, as the following argument shows. When lumping a Markov chain, the spectrum of the reduced dynamics is always a subset of the original spectrum~\cite{Barr}. For a large transition matrix, size $N \times N$, with uncorrelated random transition probabilities, the spectrum is typically concentrated to a disk with radius $\sim 1/(2  \sqrt{N})$, except for the Perron-Frobenius eigenvalue. For a block stochastic matrix the eigenvalues of the lumped Markov chain $\widetilde{P}$ are typically concentrated to a disk with radius $\sim 1/(2 \sqrt{K})$ where $K$ is the number of states in the lumped chain. If $K \ll N$ then it should be expected that the eigenvalues associated with the lumped process separate from the rest of the spectrum. An example of this can be seen in Fig.~\ref{fig:spectrum_block_stoch}. However, it should be noted that the separation is only a typical behavior, it is not necessary for the Markov chain to be lumpable (this seems to be incorrectly stated in~\cite{Shi01arandom,weinan}, see~\cite{jacobi} for details).

\begin{figure} 
	\psfrag{err}[][]{\scriptsize{$\Delta$}}
	\psfrag{eps}[][]{\scriptsize{$\epsilon$}}
	\begin{center}
		\epsfig{file=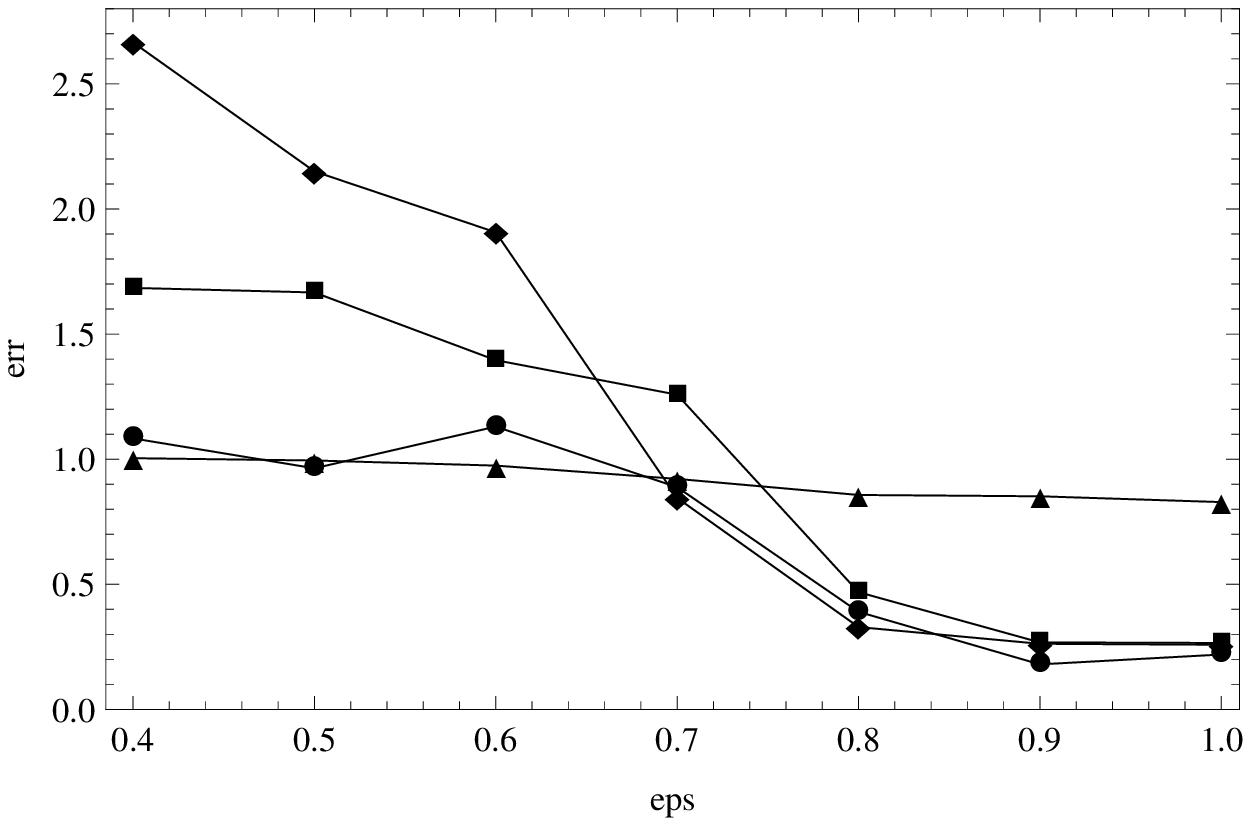,width=6cm} 
\end{center}
\caption{The average result of inferring lumping of states of a block stochastic transition matrices defined in (\ref{eq:block_stoch}). The deviation from fulfilling the lumpability condition, defined in (\ref{eq:error}), is normalized against the deviation produced by the lumping used when producing the matrix (i.e. smaller values implies better results). For each $\epsilon$ value $100$ independent realizations of $200 \times 200$ matrices on the form (\ref{eq:block_stoch}) was used to calculate the average performance. The result for the $Q$ method is displayed as $\bullet$, the result for using eigenvectors of $P$ as  $\filleddiamond$, the singular vectors from SVD \cite{friszache} as $\filledmedsquare$, and results from the clustering method suggested in~\cite{Lafon} (see Sec.~\ref{sec:lumping} for details) as $\filledmedtriangleup$. For moderate $\epsilon$ the $Q$ method and the clustering performs equally superior to the other methods, while for larger $\epsilon$ all methods except the clustering technique show approximately equal performance.}
	\label{fig:results_block_stoch}
\end{figure}	

If the spectrum does not show any separated eigenvalues that indicate the best choice of  $\lambda$ in $Q ( \lambda )$ the implementation of the $Q$ method is less straight forward when searching for general lumping. This is of course also the case for other spectral methods if the eigenvalues involved in the lumping does not separate from the rest of the spectrum. Perhaps the easiest way is to choose a set of $\{ \lambda _i \} _{i=1}^K$ randomly in the disk of radius $1$ in the complex plane, use the eigenvectors of $Q ( \lambda _i )$, $i=1, \dots , K$, corresponding to the smallest eigenvalue, cluster the elements in the same way as for the meta stable states,  and check how well the result satisfies the lumpability criterion (\ref{eq:lumping}). The procedure must be repeated a few times to find the configuration with the most satisfying result. It is possible to design more sophisticated methods by re-using the $\lambda$'s that seem to produce good results. However, we use the simplest possible approach choosing between $2$ and $5$ (guided by the separation in the spectrum) $\lambda$ values randomly in the complex plane and repeating $10$ times. The results are shown in Fig.~\ref{fig:results_block_stoch} in comparison with other methods. In these numerical test we use the deviation from 	fulfilling the lumpability condition (\ref{eq:commuting}) %
\begin{equation}
	\Delta = \| P \Pi - \Pi \widetilde{P} \|_2 
	\label{eq:error}
\end{equation}
as a measure of how well the different methods perform. For each value of $\epsilon$ test where performed with $100$ matrices generated on the form $P = (1 - \epsilon ) B + \epsilon A$, where $B$ was constructed according to (\ref{eq:block_stoch}) and $B$ was a random transition matrix. The numerical tests indicate that the $Q$ method is more stable than using $P$ directly or the SVD method. The clustering method performs almost as well as the $Q$ method.

From a computational perspective the $Q$ method is, in the case of general lumping, significantly slower than the other spectral methods. The reason is that several random choices of $\lambda$'s must be tried and in addition $Q ( \lambda )$ is a complex matrix if $\lambda$ is complex. Neither of these complications occur when searching for meta stable states since then we know beforehand that $\lambda =1$ is a good choice. In the implementation used to produce the results in Fig.~\ref{fig:results_block_stoch} the $Q$ method is approximately $15$ times slower than using $P$ directly or the SVD method. On the other hand the results are also better. The slowdown scales proportional to the number of $\lambda$-setups we need to try. A more sophisticated selection procedure for choosing the regions where the sub-dominant eigenvectors of $Q (\lambda )$ show a clear signal would probably increase the efficiency of the algorithm. 

\section{Conclusions}

We have introduced a new spectral method for identifying lumping in large Markov chains, with the identification of meta stable states as an important special case. The key element of the method is to define a family of self-adjoint matrices from the transition matrix. The eigenvectors of the self-adjoint matrices are, as opposed to those of the transition matrix itself, stable to perturbations or noisy estimation of the transition probabilities. The robustness of the method is tested and compared to the results from previous methods, including a direct clustering method introduced in ~\cite{Lafon} and the recently suggested SVD based method introduced in~\cite{friszache}. The $Q$ method is shown to be more robust than previous techniques.

We mentioned in the introduction that the method presented here can be used to reduce networks by aggregating nodes. The examples in this paper are however focused completely on lumping of Markov chains. The relation between lumping of Markov chain and reduction of complex networks was recently discussed in~\cite{weinan}. The authors define a diffusion process on the network using the standard method of the graph Laplacian. It should be noted however that reduction of networks can be defined with respect to other types of dynamics than diffusion. Straight forward jump processes are, for example, defined directly by multiplication of the adjacency matrix. Since the lumpability condition considered in this paper applies to general linear processes, not only Markov processes with stochastic transition matrix, the methods introduced can be used to reduce networks by aggregation with respect to different criteria depending on which dynamic process we are considering on the network. For the $Q$ method to be different from using the eigenvectors of the transition matrix directly, the graph must be directed.

\bibliography{articles}

\end{document}